\documentclass[12pt]{article}
\usepackage{amssymb,latexsym,amsthm,pstricks,amsmath}
\usepackage{hyperref}

\newtheorem{thm}{Theorem}[section]
\newtheorem{prop}[thm]{Proposition}
\newtheorem{lem}[thm]{Lemma}
\newtheorem{conj}[thm]{Conjecture}

\newcommand{\bth}{\begin{thm}}
\renewcommand{\eth}{\end{thm}}
\newcommand{\bcon}{\begin{conj}}
\newcommand{\econ}{\end{conj}}
\newcommand{\ble}{\begin{lem}}
\newcommand{\ele}{\end{lem}}
\newcommand{\bpr}{\begin{prop}}
\newcommand{\epr}{\end{prop}}
\newcommand{\beq}{\begin{equation}}
\newcommand{\eeq}{\end{equation}}
\newcommand{\bprf}{\begin{proof}}
\newcommand{\eprf}{\end{proof}}

\newcommand{\bbR}{{\mathbb R}}
\newcommand{\cL}{{\cal L}}
\newcommand{\cR}{{\cal R}}
\newcommand{\cH}{{\cal H}}
\newcommand{\ree}[1]{(\ref{#1})}
\newcommand{\hqed}{\hfill \qed}
\newcommand{\eqed}[1]{$\textcolor{white}{\qed}\hfill{\displaystyle#1}\hfill\qed$}
\newcommand{\fl}[1]{\lfloor #1 \rfloor}
\newcommand{\hs}[1]{\hspace{#1}}
\newcommand{\vs}[1]{\vspace{#1}}
\newcommand{\case}[4]{\left\{\begin{array}{ll}#1&\mbox{#2}\\#3&\mbox{#4}\end{array}\right.}
\newcommand{\la}{\lambda}
\newcommand{\bx}{{\bf x}}
\newcommand{\La}{\Lambda}
\newcommand{\mint}{\mathop{\rm mint}\nolimits}
\newcommand{\gaus}[2]{\genfrac{[}{]}{0pt}{}{#1}{#2}}
\newcommand{\qua}[2]{\genfrac{\langle}{\rangle}{0pt}{}{#1}{#2}}

\begin{document}

\title{Infinite log-concavity: developments and conjectures
}
\author{Peter R. W. McNamara\\[-5pt]
\small Department of Mathematics, Bucknell University,\\[-5pt]
\small Lewisburg, PA 17837, USA, \texttt{peter.mcnamara@bucknell.edu}\\
and\\
Bruce E. Sagan\\[-5pt]
\small Department of Mathematics, Michigan State University,\\[-5pt]
\small East Lansing, MI 48824-1027, USA, \texttt{sagan@math.msu.edu}
}

\date{\today\\[10pt]
	\begin{flushleft}
	                                       \\[5pt]
	\end{flushleft}}

\maketitle

\begin{abstract}
Given a sequence $(a_k)=a_0,a_1,a_2,\ldots$ of real numbers,
define a new sequence $\cL(a_k)=(b_k)$
where $b_k=a_k^2-a_{k-1}a_{k+1}$.  So $(a_k)$ is log-concave if and
only if $(b_k)$ is a nonnegative sequence.  Call $(a_k)$ {\it
infinitely log-concave\/} if $\cL^i(a_k)$ is nonnegative for all
$i\ge1$.  Boros and Moll~\cite{bm:ii} conjectured that the rows of
Pascal's triangle are infinitely log-concave.  Using a computer and a
stronger version of log-concavity, we prove their conjecture for the
$n$th row for all $n\le 1450$.  We also use our methods to give a
simple proof of a 
recent result of Uminsky and Yeats~\cite{uy:uri} about regions of infinite
log-concavity.  We investigate related questions about the
columns of Pascal's triangle, $q$-analogues, symmetric
functions, real-rooted polynomials, and
Toeplitz matrices.  In addition, we offer several conjectures.
\end{abstract}

\section{Introduction}

Let 
$$
(a_k)=(a_k)_{k\ge0}=a_0,a_1,a_2,\ldots
$$ 
be a sequence of real numbers.  It will
be convenient to extend the sequence to negative indices by letting
$a_k=0$ for $k<0$. Also, if $(a_k)=a_0,a_1,\ldots,a_n$ is a finite
sequence then we let $a_k=0$ for $k>n$.

Define the {\it $\cL$-operator\/} on sequences to be $\cL(a_k)=(b_k)$
where $b_k=a_k^2-a_{k-1}a_{k+1}$.  
Call a sequence {\it $i$-fold log-concave\/} if $\cL^i(a_k)$
is a nonnegative sequence.  So log-concavity in the ordinary sense is
$1$-fold log-concavity.
Log-concave sequences arise in many areas of algebra, combinatorics,
and geometry.  See the survey articles of Stanley~\cite{sta:lus} and
Brenti~\cite{bre:lus} for more information.

Boros and Moll~\cite[page 157]{bm:ii} defined $(a_k)$ to be 
{\it  infinitely log-concave\/} if it is $i$-fold log-concave  for all
$i\ge1$.  They introduced this definition in conjunction with the
study of a specialization of the Jacobi polynomials whose coefficient
sequence they conjectured to be infinitely log-concave.   
Kauers and Paule~\cite{kp:cpm} used a computer algebra package to prove
this conjecture for ordinary log-concavity.  
Since the coefficients of these polynomials can be expressed in terms of
binomial coefficients, Boros and Moll also made the statement:
\begin{center}
``Prove that the binomial coefficients are $\infty$-logconcave.''
\end{center}
We will take this to be a conjecture that the rows of Pascal's triangle are infinitely log-concave, although we will later discuss the columns and other lines. 
When
given a function of more than one variable, we will 
subscript the $\cL$-operator by the parameter which is varying to form
the sequence.  
So $\cL_k\binom{n}{k}$ would refer to the operator acting on the
sequence $\binom{n}{k}_{k\ge0}$.  
Note that we drop the sequence parentheses
for sequences of binomial coefficients to improve readability.
We now restate
the Boros-Moll conjecture formally.
\bcon
\label{bm:con}
The sequence $\binom{n}{k}_{k\ge0}$ is  infinitely log-concave for all
$n\ge0$.
\econ

In the next section, we use a strengthened version of log-concavity
and computer calculations to verify Conjecture~\ref{bm:con} for all
$n\le 1450$.  Uminsky and Yeats~\cite{uy:uri} set up a correspondence
 between certain symmetric
sequences and points in $\bbR^m$.  They then described an infinite region 
$\cR\subset\bbR^m$ bounded by hypersurfaces and such that
each sequence corresponding to a point of $\cR$ 
is infinitely log-concave.  In Section~\ref{ril}, we show how
our methods can be used to give a simple derivation of one of their
main theorems.  We investigate infinite log-concavity of the columns
and other lines of Pascal's 
triangle in Section~\ref{cpt}.  
Section~\ref{gps} is devoted to two $q$-analogues of the binomial
coefficients.  For the Gaussian polynomials, we show that certain
analogues of some infinite log-concavity conjectures 
are false while others appear to be true.  In contrast, our second
$q$-analogue seems to retain all the log-concavity properties of the
binomial coefficients.   In Section~\ref{sf}, after showing why the
sequence $(h_k)_{k\ge0}$ of complete homogeneous symmetric is an
appropriate analogue of sequences of binomial coefficients, we explore
its log-concavity properties.  We end with a section of related
results and questions about real-rooted polynomials and Toeplitz matrices.

While one purpose of this article is to present our
  results, we have written it with two more targets in mind.
  The first is to convince our audience that infinite log-concavity is
  a fundamental concept.  We hope that its definition as a natural
  extension of traditional log-concavity helps to make this case. 
Our second aspiration is to attract others to work on the subject; to
  that end, we have presented several open problems.  These
  conjectures each represent fundamental questions in the area, so even
  solutions of special cases may be interesting. 

\section{Rows of Pascal's triangle}
\label{rpt}

One of the difficulties with proving the Boros-Moll conjecture is that
log-concavity is not preserved by the $\cL$-operator.  For example,
the sequence 
$4,5,4$ is log-concave but $\cL(4,5,4)=16,9,16$ is not.  So we
will seek a condition stronger than log-concavity which is preserved
by $\cL$.  Given $r\in\bbR$, we say that a sequence $(a_k)$ is
{\it $r$-factor  log-concave\/} if 
\beq
\label{r}
a_k^2\ge r a_{k-1}a_{k+1}
\eeq
for all $k$.  Clearly this implies log-concavity if $r\ge1$.

We seek an $r>1$ such that $(a_k)$ being $r$-factor log-concave implies
that $(b_k)=\cL(a_k)$ is as well.  Assume the original sequence is
nonnegative.  Then expanding $r b_{k-1} b_{k+1}\le b_k^2$ in terms of
the $a_k$ and rearranging the summands, we see that this is equivalent
to proving
$$
(r-1) a_{k-1}^2 a_{k+1}^2 + 2 a_{k-1} a_k^2 a_{k+1}
\le a_k^4 + r  a_{k-2} a_k (a_{k+1}^2-a_k a_{k+2} ) + r a_{k-1}^2 a_k a_{k+2}.
$$

By our assumptions, the two expressions with factors of $r$ on the right are
nonnegative, so it suffices to prove the inequality obtained when
these are dropped.  Applying~\ree{r} to the left-hand side gives
$$
(r-1) a_{k-1}^2 a_{k+1}^2 + 2 a_{k-1} a_k^2 a_{k+1}
\le  \frac{r-1}{r^2} a_k^4 + \frac{2}{r} a_k^4.
$$
So we will be done if
$$
 \frac{r-1}{r^2}+ \frac{2}{r} =1.
$$
Finding the root $r_0>1$ of the corresponding quadratic equation
finishes the proof of the first assertion of the following lemma, while the second assertion follows easily from the first.

\ble
\label{rfactor:le}
Let 
$(a_k)$ be a nonnegative sequence and let $r_0=(3+\sqrt{5})/2$.  
Then $(a_k)$ being $r_0$-factor log-concave 
implies that $\cL(a_k)$ is too. So in this case $(a_k)$ is infinitely
log-concave.\hqed
\ele

Now to prove that any row of Pascal's triangle is infinitely
log-concave, one merely lets a computer find $\cL_k^i\binom{n}{k}$ for
$i$ up to some bound $I$.
If these sequences are all nonnegative and
$\cL_k^I\binom{n}{k}$ is $r_0$-factor log-concave, then the
previous lemma shows that this row is infinitely log-concave.  Using
this technique, we have obtained the following theorem.
\bth
\label{row:th}
The sequence $\binom{n}{k}_{k\ge0}$ is infinitely log-concave for
all $n\le 1450$.\hqed
\eth

We note that the necessary value of $I$ increases slowly with increasing $n$.
As an example, when $n=100$, our technique works with $I=5$, while for
$n=1000$, we need $I=8$.

Of course, the method developed in this section can be applied to any
sequence such that $\cL^i(a_k)$ is $r_0$-factor log-concave for some
$i$.  In particular, it is interesting to try it on the original
sequence which motivated Boros and Moll~\cite{bm:ii} to define infinite
log-concavity.  They were studying the polynomial
\beq
\label{Pm}
P_m(x) = \sum_{\ell=0}^m d_\ell(m) x^\ell
\eeq
where
$$
d_\ell(m) = \sum_{j=\ell}^m 2^{j-2m} 
\binom{2m-2j}{m-j}\binom{m+j}{m}\binom{j}{\ell}.
$$
Kauers~[private communication] has used our method to verify
infinite log-concavity of the sequence $(d_\ell(m))_{\ell\ge0}$ for
$m\le 129$.  For such values of $m$, $\cL^5_\ell$ applied to the sequence is 
$r_0$-factor log-concave.

\section{A region of infinite log-concavity}
\label{ril}

Uminsky and Yeats~\cite{uy:uri}  took a different approach to the
Boros-Moll Conjecture as described in the Introduction.  Since they
were motivated by the rows of Pascal's triangle, they only considered
real sequences $a_0,a_1,\ldots,a_n$ which are symmetric (in that
$a_k=a_{n-k}$ for all $k$) and satisfy $a_0=a_n=1$.  Each such sequence
corresponds to a point $(a_1,\ldots,a_m)\in\bbR^m$ where $m=\fl{n/2}$.

Their region, $\cR$, whose points all correspond to infinitely log-concave
sequences, is bounded by $m$ parametrically defined hypersurfaces.  The
parameters are $x$ and $d_1,d_2,\ldots,d_m$  and it will be convenient
to have the notation
$$
s_k=\sum_{i=1}^k d_i.
$$
We will also need $r_1=(1+\sqrt{5})/2$.  Note that $r_1^2=r_0$.
The $k$th hypersurface, $1\le k< m$, is defined as
$$
\begin{array}{l}
\cH_k=
\{(x^{s_1},\ldots,x^{s_{k-1}},r_1x^{s_k},x^{s_{k+1}+d_k-d_{k+1}},\ldots,x^{s_m+d_k-d_{k+1}}):
\\[5pt]
\hs{100pt} x\ge1, \quad 1=d_1>\cdots> d_k>d_{k+2}>\cdots>d_m>0\},
\end{array}
$$
while
$$
\cH_m=
\{(x^{s_1},\ldots,x^{s_{m-1}},cx^{s_{m-1}}):\
x\ge1, \quad 1=d_1>\cdots > d_{m-1}>0\},
$$
where
$$
c=\case{r_1}{if $n=2m$,}{2}{if $n=2m+1$.}
$$

Let us say that the \emph{correct side} of $\cH_k$ for $1\le k\le m$ consists of those points in $\bbR^m$ that can be obtained from a point on $\cH_k$ by increasing the $k$th coordinate.
Then let $\cR$ be the region of all points in $\bbR^m$ having
increasing coordinates and lying on the correct side of $\cH_k$ for all $k$.
We will show how our method of the previous section can be used to give
a simple proof  of one of  Uminsky and Yeats' main theorems.  But first we
need a modified version of Lemma~\ref{rfactor:le} to take care of the
case when $n=2m+1$.

\ble
\label{rfactor:mod}
Let $a_0,a_1,\ldots, a_{2m+1}$ be a symmetric, nonnegative sequence such that
\begin{enumerate}
\item[(i)]  $a_k^2 \ge r_0 a_{k-1} a_{k+1}$ for $ k< m$, and
\item[(ii)]  $a_m\ge 2 a_{m-1}$.
\end{enumerate}
Then $\cL(a_k)$ has the same properties, which implies that $(a_k)$ is
infinitely log-concave.
\ele
\bprf
Clearly $\cL(a_k)$ is still symmetric.  To show that the other
two properties persist, note that in demonstrating
Lemma~\ref{rfactor:le} we actually proved more.  In particular, we
showed that if equation~\ree{r} holds at index $k$ of the 
sequence $(a_k)$ (with $r=r_0$), then it also holds at index $k$ of the
sequence $\cL(a_k)$ provided that the original sequence is log-concave.
Note that the assumptions of the current lemma imply log-concavity of
$(a_k)$:  This is clear at indices $k\neq m,m+1$ because of 
condition~(i).  Also, using symmetry 
and multiplying condition (ii) by $a_m$ gives
$a_m^2 \ge 2 a_{m-1} a_{m} = 2 a_{m-1} a_{m+1}$ (and symmetrically for
$k=m+1$).  

So now we know that condition~(i) is also true for $\cL(a_k)$.  As for
condition~(ii), using symmetry we see that we need to prove
$$
a_m^2-a_{m-1} a_m\ge 2\left( a_{m-1}^2-a_{m-2} a_m \right).
$$
Rearranging terms and dropping one of them shows that it suffices to
demonstrate
$$
2 a_{m-1}^2 + a_{m-1} a_m \le a_m^2.
$$
But this is true because of~(ii), and we are done.
\eprf

\bth[\cite{uy:uri}]
Any sequence corresponding to a point of $\cR$ is infinitely log-concave.
\eth
\bprf
It suffices to show that the sequence satisfies the hypotheses of 
Lemma~\ref{rfactor:le} when $n=2m$, or Lemma~\ref{rfactor:mod} when $n=2m+1$.

Suppose first that
$k<m$.  Being on the correct side of $\cH_k$ 
is equivalent to there being values of the parameters such that
$$
a_k^2
\hs{5pt}\ge\hs{5pt} 
(r_1 x^{s_k})^2
\hs{5pt}=\hs{5pt}
r_1^2 x^{(s_{k-1}+d_k)+(s_{k+1}-d_{k+1})}
\hs{5pt}=\hs{5pt}
r_0 a_{k-1} a_{k+1}.
$$
Thus we have the necessary inequalities for this range of $k$.

If $k=m$ then we can use an argument as in the previous paragraph if $n=2m$.
If
$n=2m+1$, then being on the correct side of $\cH_m$ 
is equivalent to
$$
a_m\ge 2 x^{s_{m-1}}=2 a_{m-1}.
$$
This is precisely condition~(ii) of Lemma~\ref{rfactor:mod}, which finishes
the proof.
\eprf

\section{Columns and other lines of Pascal's triangle} 
\label{cpt}

While we have treated Boros and Moll's statement about the infinite
log-concavity of the binomial coefficients to be a statement about the
rows of Pascal's triangle, their wording also suggests an examination
of the columns.
\bcon
\label{col:con}
The sequence $\binom{n}{k}_{n\ge k}$ is  infinitely log-concave for
all fixed $k\ge0$. 
\econ

We will give two pieces of evidence for this conjecture.  One is a demonstration
that various columns corresponding to small values 
of $k$ are infinitely log-concave.
 Another is a proof that $\cL_n^i\binom{n}{k}$ is
nonnegative for certain values of $i$ and all $k$.

\bpr
\label{small:k}
The sequence $\binom{n}{k}_{n\ge k}$ is  infinitely log-concave for
$0\le k \le 2$.
\epr
\bprf  When $k=0$ we have, for all $i \geq 1$, 
$$
\cL_n^i\binom{n}{0}=(1,0,0,0,\ldots).
$$
For $k=1$ we obtain
$$
\cL_n\binom{n}{1} = (1,1,1,\dots)
$$
so infinite log-concavity follows from the $k=0$ case.
The sequence when $k=2$ is a fixed point of
the $\cL$-operator, again implying infinite log-concavity.   
\eprf

In what follows, we use the notation $L(a_k)$ for the $k$th
element of the sequence $\cL(a_k)$, and similarly for $L_k$ and $L_n$.

\bpr
\label{small:i}
The sequence $\cL_n^i\binom{n}{k}$ is nonnegative for all $k$ and for
$0\le i\le 4$.
\epr
\bprf
By the previous proposition, we only need to check $k\ge3$.
Using the expression for a binomial coefficient in terms of
factorials, it is easy to derive the following expressions:
$$
L_n\binom{n}{k} = 
\frac{1}{n}\binom{n}{k}\binom{n}{k-1}
$$
and
$$
L_n^2\binom{n}{k} = 
\frac{2}{n^2 (n-1)} 
\binom{n}{k}^2 \binom{n}{k-1}\binom{n}{k-2}.
$$

With a little more work, one can show that $L_n^3\binom{n}{k}$ can
be expressed as a product of nonnegative factors times the polynomial
$$
(4k-6)n^2-(4k^2-10k+6)n-k^2.
$$ 
To show that this is nonnegative, we write $n=k+m$ for $m\ge0$ to get
$$
(4k-6)m^2+(4k^2-2k-6)m+(3k^2-6k).
$$
But the coefficients of the powers of $m$ are all positive for $k\ge3$,
so we are done with the case $i=3$.

When $i=4$, we follow the same procedure, only now the polynomial in
$m$ has coefficients which are polynomials in $k$ up to degree $7$.
For example, the coefficient of $m^3$ is
$$
528 k^7- 8 k^6 - 11,248 k^5 + 25,360 k^4 - 5,888 k^3 - 24,296 k^2
+16,080 k - 1,584. 
$$
To make sure this is nonnegative for integral $k\ge3$, one rewrites
the polynomial as
$$
(528 k^2- 8 k - 11,248) k^5 + (25,360 k^2 - 5,888 k - 24,296) k^2
+(16,080 k - 1,584), 
$$
finds the smallest $k$ such that each of the factors in parentheses is
nonnegative from this value on, and then checks any remaining $k$ by
direct substitution. 
\eprf

Kauers and Paule~\cite{kp:cpm} proved that the rows of Pascal's
triangle are $i$-fold log-concave for $i\le5$.
Kauers~[private communication] has used their techniques to confirm 
Proposition~\ref{small:i} and to also check the case $i=5$ for the
columns.  For the latter case, Kauers used a computer to determine 
\beq 
\label{kau}
\frac{(\cL_n^5 \binom{n}{k})}{\binom{n}{k}^{2^5}}
\eeq
explicitly, which 
is just a rational function in $n$ and $k$.  He then showed that  
\ree{kau} is nonnegative by means of cylindrical algebraic decomposition.
We refer the interested reader to~\cite{kp:cpm} and the references therein for more
information on such techniques.

More generally, we can look at an arbitrary line in Pascal's triangle,
i.e.,  consider the sequence $\binom{n+mu}{k+mv}_{m\ge0}$.  
The unimodality and (1-fold) log-concavity of such sequences has been investigated in 
\cite{bbs:ucs, sw:upp, tz:usb, tz:usb2}.
We do not require that $u$ and $v$ be coprime, so such sequences
need not contain all of the binomial coefficients in which a geometric
line would intersect Pascal's triangle, e.g., a sequence such as 
$\binom{n}{0}, \binom{n}{2},\binom{n}{4},\ldots$ would be included.
By letting $u<0$, one can get a finite truncation of a column.  For
example, if $n=5$, $k=3$, $u=-1$, and $v=0$ then we get the sequence
$$
\binom{5}{3}, \binom{4}{3}, \binom{3}{3}
$$ 
which is not even $2$-fold log-concave.   So we will only consider $u\ge0$.
Also
$$
\binom{n+mu}{k+mv} = \binom{n+mu}{n-k+m(u-v)}
$$  
so we can also assume $v\ge0$.

We offer the following conjecture, which includes 
Conjecture~\ref{bm:con} as a special case.
\bcon
\label{all:con}
Suppose that $u$ and $v$ are distinct nonnegative integers.  Then
$\binom{n+mu}{mv}_{m\ge0}$ is infinitely log-concave for all $n \geq 0$ if
and only if $u<v$ or $v=0$. 
\econ

We first give a quick proof of the ``only if'' direction.  Supposing
that $u>v \ge1$, we consider the sequence 
$$
\binom{0}{0}, \binom{u}{v}, \binom{2u}{2v}, \ldots
$$
obtained when $n=0$.
We claim that this sequence is not even log-concave and that
log-concavity fails at the second term.  Indeed, the fact that
$\binom{u}{v}^2 < \binom{2u}{2v}$ follows immediately from the
identity
\[
\binom{u}{0}\binom{u}{2v} + \binom{u}{1}\binom{u}{2v-1} + \cdots + 
\binom{u}{v}\binom{u}{v} + \cdots + \binom{u}{2v}\binom{u}{0} = \binom{2u}{2v},
\]
which is a special case of Vandermonde's Convolution.

The proof just given shows that subsequences of the columns of
Pascal's triangle are the only {\it infinite} sequences of the form
$\binom{n+mu}{mv}_{m\ge0}$ that can possibly be infinitely log-concave.   
We also note that the previous conjecture says nothing about what
happens on the diagonal $u=v$.  Of course, the case $u=v=1$ is
Conjecture~\ref{col:con}.  For other diagonal values, the evidence is
conflicting. 
One can show by computer that $\binom{n+mu}{mu}_{m\ge0}$ is not $4$-fold log-concave
for $n=2$ and any
$2\le u\le500$.  However, this is the only {\it known} value of $n$
for which $\binom{n+mu}{mu}_{m\ge0}$ is not an infinitely log-concave
sequence for some $u\ge1$. 

We conclude this section by offering considerable computational
evidence in favor of the ``if'' direction of Conjecture~\ref{all:con}.
Theorem~\ref{row:th} provides such evidence when $u=0$ and $v=1$.  Since all
other
sequences with $u<v$ have a finite number of nonzero entries, we can
use the $r_0$-factor log-concavity technique for these sequences as well.
For all $n\le 500$, $2\le v\le20$ and $0\le u<v$, we have checked that
$\binom{n+mu}{mv}_{m\ge0}$ is infinitely log-concave.

\section{$q$-analogues}
\label{gps}

This section will be devoted to discussing two $q$-analogues of
binomial coefficients.   For the Gaussian polynomials, we will see
that the corresponding generalization of Conjecture~\ref{bm:con} is
false, and we show one exact reason why it fails.  In contrast, the
corresponding generalization of Conjecture~\ref{col:con} appears to be
true.   This shows how delicate these conjectures are and may in part
explain why they seem to be difficult to prove.  After introducing our
second $q$-analogue, we
conjecture that the corresponding generalizations of
Conjectures~\ref{bm:con}, \ref{col:con} and \ref{all:con} are all
true.  This second $q$-analogue arises in the study of quantum groups; see,
for example, the books of Jantzen~\cite{jan:lqg} and Majid~\cite{maj:qgp}.

Let $q$ be a variable and consider a polynomial $f(q)\in\bbR[q]$.
Call $f(q)$ \linebreak $q$-{\it nonnegative\/}  if all the coefficients of $f(q)$ are
nonnegative.  Apply
the $\cL$-operator to sequences of polynomials $(f_k(q))$ in the
obvious way.  Call such a sequence {\it $q$-log-concave\/} if
$\cL(f_k(q))$ is a sequence of $q$-nonnegative polynomials, with
$i$-fold $q$-log-concavity and infinite $q$-log-concavity defined
similarly.  

We will be particularly interested in the Gaussian polynomials.  The
standard {\it $q$-analogue of the nonnegative integer $n$\/} is
$$
[n]=[n]_q=\frac{1-q^n}{1-q} = 1+q+q^2+\cdots+q^{n-1}.
$$
Then, for $0\le k\le n$, the {\it Gaussian polynomials\/} or 
{\it $q$-binomial coefficients\/} are defined as
$$
\gaus{n}{k}=\gaus{n}{k}_q = \frac{[n]_q!}{[k]_q! [n-k]_q!}
$$
where $[n]_q!=[1]_q [2]_q\cdots [n]_q$.   For more information,
including proofs of the assertions made in the next paragraph, see the
book of Andrews~\cite{and:tp}.

Clearly substituting $q=1$ gives $\genfrac{[}{]}{0pt}{}{n}{k}_1=\binom{n}{k}$.
Also, it is well known that the Gaussian polynomials
are indeed $q$-nonnegative polynomials.  In fact, they have
various combinatorial interpretations, one of which we will need.
An (integer) {\it partition of $n$\/} is a weakly decreasing
positive integer sequence $\la=(\la_1,\la_2,\ldots,\la_\ell)$ such
that $|\la|\stackrel{\rm def}{=}\sum_i \la_i=n$.  The $\la_i$ are
called {\it parts\/}.  For notational convenience, if a part $k$ is repeated $r$ times in a
partition $\la$ then we will denote this by writing $k^r$ in the
sequence for $\la$.
We say that $\la$ {\it fits inside an $s\times t$ box} if
$\la_1\le t$ and  $\ell\le s$.  Denote the set of all such partitions
by $P(s,t)$. 
It is well known, and easy to prove by induction on $n$, that
\beq
\label{gau}
\gaus{n}{k}=\sum_{\la\in P(n-k,k)} q^{|\la|}.
\eeq

We are almost ready to prove that the sequence
$\left(\gaus{n}{k}\right)_{k\ge0}$ is not infinitely $q$-log-concave.
In fact, we will show it is not even $2$-fold $q$-log-concave.
First we need a lemma.  In it, we use $\mint f(q)$ to denote the
nonzero term of least degree in $f(q)$.
\ble
\label{mint}
Let $L_k\left(\gaus{n}{k}\right)=B_k(q)$.  Then for $k\le n/2$,
$$
\mint B_k(q)=
\case{q^k}{if $k<n/2$,}{2q^k}{if $k=n/2$.}
$$
\ele
\bprf
Since $B_k(q)=\gaus{n}{k}^2-\gaus{n}{k-1}\gaus{n}{k+1}$ it suffices to
prove, in view of~\ree{gau}, the following two statements.  If $i\le k$ and 
$$
(\la,\mu)\in P(n-k+1,k-1)\times P(n-k-1,k+1)
$$ 
with $|\la|+|\mu|=i$, then $(\la,\mu)\in P(n-k,k)^2$.  Furthermore,
the number of elements in 
$P(n-k,k)^2- P(n-k+1,k-1)\times P(n-k-1,k+1)$ is 0 or 1 or 2 depending
on whether $i<k$ or $i=k<n/2$ or $i=k=n/2$, respectively.

The first statement is an easy consequence of 
$|\la|+|\mu|=i\le k\le n-k$.  A similar argument works for the $i<k$
case of the second statement.  If $i=k$ then the pair $((k),\emptyset)$ is in the
difference and if $i=k=n/2$ then the pair $(\emptyset,(1^k))$ is as well.
\eprf

\bpr
\label{gau:pr}
Let $L_k^2\left(\gaus{n}{k}\right)=C_k(q)$.  Then for
$n\ge2$ and $k=\fl{n/2}$ we have
$$
\mint C_k(q)=-q^{n-2}.
$$
Consequently, $\left(\gaus{n}{k}\right)_{k\ge0}$ is not $2$-fold $q$-log-concave.
\epr
\bprf
The proofs for $n$ even and odd are similar, so we will only do the
former.  So suppose $n=2k$ and consider 
$$
C_k(q)=B_k(q)^2-B_{k-1}(q)B_{k+1}(q)=B_k(q)^2-B_{k-1}(q)^2.
$$
By the previous lemma
$\mint B_k(q)^2=4q^{2k}$ and $\mint B_{k-1}(q)^2= q^{2k-2}$.  Thus
$\mint C_k(q) = -q^{2k-2}=-q^{n-2}$ as desired.
\eprf

After what we have just proved, it may seem surprising that the
following conjecture, which is a $q$-analogue of
Conjecture~\ref{col:con}, does seem to hold. 
\bcon
\label{qcol:con}
The sequence $\left(\gaus{n}{k}\right)_{n\ge k}$ is  infinitely
$q$-log-concave for
all fixed $k\ge0$. 
\econ

As evidence, we will prove a $q$-analogue of Proposition~\ref{small:k}
and comment on Proposition~\ref{small:i} in this setting.
\bpr
\label{qsmall:k}
The sequence $\left(\gaus{n}{k}\right)_{n\ge k}$ is  infinitely
$q$-log-concave for
$0\le k\le 2$.
\epr
\bprf
When $k=0$ one has the same sequence as when $q=1$.  

When $k=1$ we claim that
$$
\cL\left(\gaus{n}{1}\right) = (1,q,q^2,q^3,\ldots).
$$
Indeed, 

\begin{eqnarray*}
[n]^2-[n-1][n+1]
&=& \frac{(1-q^n)^2-(1-q^{n-1})(1-q^{n+1})}{(1-q)^2}\\
&=& \frac{q^{n-1}-2q^n+q^{n+1}}{(1-q)^2}\\
&=& q^{n-1}
\end{eqnarray*}
(and recall that the sequence starts at $n=1$).  It follows that
$$
\cL^i\left(\gaus{n}{1}\right) = (1,0,0,0,\ldots)
$$
for $i\ge2$.

For $k=2$, the manipulations are much like those in the previous paragraph.
Using induction on $i$, we obtain
$$
L^i\left(\gaus{n}{2}\right) = q^{(2^i-1)(n-2)}\gaus{n}{2}
$$
for $i\ge0$.  This completes the proof of the last case of the
proposition.  
\eprf

If we now consider arbitrary $k$ it is not hard to show, using
algebraic manipulations like those in the proof just given,  
that
\beq
\label{lq}
L_n \left(\gaus{n}{k}\right)=
\frac{q^{n-k}}{[n]}\gaus{n}{k}\gaus{n}{k-1}.
\eeq
These are, up to a power of $q$, the $q$-Narayana numbers.  They were
introduced by F\"urlinger and Hofbauer~\cite{fh:qcn} and are contained
in a specialization of a result of MacMahon~\cite[page 1429]{mac:cp}
which was stated without proof.  They were further studied by
Br\"and\'en~\cite{bra:qnn}.  As shown in the references just cited,
these polynomials are the generating functions for a number of
different families of combinatorial objects.  Thus they are $q$-nonnegative.

More computations show that
\beq
\label{llq}
L_n^2\left(\gaus{n}{k}\right)
=\frac{q^{3n-3k}[2]}{[n]^2[n-1]}\gaus{n}{k}^2\gaus{n}{k-1}\gaus{n}{k-1}.
\eeq
It is not clear that these polynomials are $q$-nonnegative, although they
must be if Conjecture~\ref{qcol:con} is true.  Furthermore, when $q=1$,
the triangle made as $n$ and $k$ vary is not in Sloane's
Encyclopedia~\cite{slo:oei} (although it has now been submitted).  We
expect that these integers and polynomials have interesting, yet to be
discovered, properties.

We conclude our discussion of the Gaussian polynomials by considering the
sequence 
\beq
\label{gs}
\left(\gaus{n+mu}{mv}\right)_{m\ge0}
\eeq
for nonnegative integers $u$ and $v$, as we did in Section~\ref{cpt} for the
binomial coefficients.
When $u>v$ the sequence has an infinite number of nonzero entries.  We
can use $\eqref{gau}$ to show that the highest degree term in
$\gaus{n+u}{v}^2 - \gaus{n+2u}{2v}$ has coefficient $-1$, so the
sequence \eqref{gs} is not even $q$-log-concave.   When $u<v$, it
seems to be the case that the sequence is not $2$-fold $q$-log-concave,
as shown for the rows in Proposition~\ref{gau:pr}.  When $u=v$, the
evidence is conflicting, reflecting the behavior of the binomial
coefficients.  Since setting $q=1$ in $\gaus{n+mu}{mu}$ yields 
$\binom{n+mu}{mu}$, we know that
$\left(\gaus{2+mu}{mu}\right)_{m\ge0}$ is not always 
$4$-fold $q$-log-concave.  It also transpires that the case $n=3$ is
not always $5$-fold $q$-log-concave.  We have not encountered other
values of $n$ that fail to yield a $q$-log-concave sequence when
$u=v$. 

While the variety of behavior of the Gaussian polynomials is
interesting, it would be desirable to have a $q$-analogue that
better reflects the behavior of the binomial coefficients.  A
$q$-analogue that arises in the study of quantum groups serves
this purpose.  Let us 
replace the previous $q$-analogue of the nonnegative integer $n$ 
with the expression 
$$
\langle n\rangle = \frac{q^n-q^{-n}}{q-q^{-1}} =
q^{1-n} + q^{3-n} + q^{5-n} + \cdots + q^{n-1}.
$$
From this, we obtain a $q$-analogue of the binomial coefficients by
proceeding as for the Gaussian polynomials: for $0\le k\le n$, we
define 
$$
\qua{n}{k} = \frac{\langle n\rangle!}{\langle k\rangle! \langle n-k\rangle!}
$$
where $\langle n\rangle!=\langle1 \rangle \langle2 \rangle\cdots \langle
n\rangle$. 

Letting $q\to1$ in $\qua{n}{k}$ gives $\binom{n}{k}$, and a
straightforward calculation shows that  
\beq
\label{quagau}
\qua{n}{k} = \frac{1}{q^{nk-k^2}}\gaus{n}{k}_{q^2}.
\eeq
So $\qua{n}{k}$ is, in general, a Laurent polynomial in $q$
with nonnegative coefficients.  Our definitions of $q$-nonnegativity and $q$-log-concavity for polynomials in $q$ extend to Laurent polynomials in the obvious way.

We offer the following generalizations of Conjectures~\ref{bm:con},
\ref{col:con} and \ref{all:con}. 

\bcon
\label{qua:con}
\
\begin{itemize} 
\item[(a)] The row sequence $\left(\qua{n}{k}\right)_{k\ge0}$ is
  infinitely $q$-log-concave for all $n\ge0$. 
\item[(b)] The column sequence $\left(\qua{n}{k}\right)_{n\ge
    k}$ is infinitely $q$-log-concave for all fixed $k\ge0$. 
\item[(c)] For all integers $0\le u<v$, the sequence
  $\left(\qua{n+mu}{mv}\right)_{m\ge0}$ is infinitely
  $q$-log-concave for all $n\ge0$. 
\end{itemize}
\econ

Several remarks are in order.  Suppose that for $f(g), g(q) \in
\bbR[q, q^{-1}]$, we say $f(q) \le g(q)$  if $g(q)-f(q)$ is $q$-nonnegative.
Then the proofs of Lemmas~\ref{rfactor:le} and
\ref{rfactor:mod} work equally well if the $a_i$'s are Laurent
polynomials and we replace the term ``log-concave'' by
``$q$-log-concave.''  Using these lemmas, we have verified
Conjecture~\ref{qua:con}(a) for all $n\le 53$.   Even though (a) is a
special case of (c), we state it separately since (a) is the
$q$-generalization of the Boros-Moll conjecture, the primary
motivation for this paper.  

As evidence for Conjecture~\ref{qua:con}(b), it is not hard to prove
the appropriate analogue of Propositions~\ref{small:k} and
\ref{qsmall:k}, i.e.\ that the sequence $\qua{n}{k}_{n\ge k}$ is
infinitely $q$-log-concave for all $0\le k\le2$.  To obtain the
expressions for $L_n\left(\qua{n}{k}\right)$ and
$L_n^2\left(\qua{n}{k}\right)$, take
equations~\eqref{lq} and \eqref{llq}, replace all square
brackets by angle brackets and replace each the terms $q^{n-k}$ and
$q^{3n-3k}$ by the number 1. 

Conjecture~\ref{qua:con}(c) has been verified for all $n\le 24$
with $v\le10$.   When $u > v$, we can use $\eqref{quagau}$ to show
that the lowest degree term in  $\qua{n+u}{v}^2 - \qua{n+2u}{2v}$ has
coefficient $-1$, so the sequence is not even $q$-log-concave. 
When $u=v$, the quantum groups analogue has exactly the same behavior
as we observed above for the Gaussian polynomials. 

\section{Symmetric functions}
\label{sf}

We now turn our attention to symmetric functions.  
We will demonstrate that the complete homogeneous symmetric functions
$(h_k)_{k\ge0}$ are a natural analogue of the rows and columns of
Pascal's triangle.  We show that the sequence
$(h_k)_{k\ge0}$ is 
$i$-fold log-concave in the appropriate sense for $i\le3$, 
but not $4$-fold log-concave.
Like the results of Section~\ref{gps}, this
result underlines the difficulties and subtleties of
Conjectures~\ref{bm:con} and~\ref{col:con}.  In particular, it shows
that any proof of Conjecture~\ref{bm:con} or Conjecture~\ref{col:con}
would need to use techniques that do not carry over to the sequence
$(h_k)_{k\ge0}$.   
For a more detailed
exposition of the background material below, we refer the reader to
the texts of Fulton~\cite{ful:yt}, Macdonald~\cite{mac:sfh},
Sagan~\cite{sag:sym} or Stanley~\cite{sta:ec2}.

Let $\bx=\{x_1,x_2,\ldots\}$ be a countably infinite set of variables.
For each $n\ge0$, the elements of the symmetric group $\mathfrak{S}_n$ act on formal
power series $f(\bx)\in\bbR[[\bx]]$ by permutation of variables (where
$x_i$ is left fixed if $i>n$).  The algebra of symmetric functions,
$\La(\bx)$, is the set of all series left fixed by all symmetric
groups and of bounded (total) degree.  

The vector space of symmetric functions homogeneous of degree $k$ has
dimension equal to the number of partitions
$\la=(\la_1,\ldots,\la_\ell)$ of $k$.  We will be 
interested in three bases for this vector space.  The 
{\it monomial symmetric function\/} corresponding to $\la$,
$m_\la=m_\la(\bx)$, is obtained by 
symmetrizing the monomial $x_1^{\la_1}\cdots x_\ell^{\la_\ell}$.  The
$k$th {\it complete homogeneous symmetric function\/}, $h_k$, is the sum of
all monomials of degree $k$.  For partitions, we then define 
$$
h_\la=h_{\la_1}\cdots h_{\la_\ell}.
$$  
Finally, the 
{\it Schur function\/} corresponding to $\la$ is
$$
s_\la= \det(h_{\la_i-i+j})_{1\le i,j \le \ell}.
$$
We remark that this determinant is a minor of the Toeplitz matrix for
the sequence $(h_k)$.  We will have more to say about Toeplitz matrices
in the next section.

Our interest will be in the sequence just mentioned
$(h_k)_{k\ge0}$.  Let $h_k(1^n)$ 
denote the integer obtained by substituting $x_1=\cdots=x_n=1$ and
$x_i=0$ for $i>n$ into $h_k=h_k(\bx)$.  Then $h_k(1^n)=\binom{n+k-1}{k}$
(the number of ways of choosing $k$ things from $n$ things with
repetition) and so the above sequence becomes a column of Pascal's
triangle.  By the same token $h_k(1^{n-k})=\binom{n-1}{k}$ and so the
sequence becomes a row.

We will now collect the results from the theory of symmetric functions
which we will need.  
Partially order partitions by {\it dominance\/} 
where $\la\le\mu$ if and only if for every $i\ge1$ we have
$\la_1+\cdots+\la_i\le\mu_1+\cdots+\mu_i$.  Also,  if  $\{b_\la\}$ is
any basis of $\La(\bx)$ and $f\in\La(\bx)$ then we let
$[b_\la]f$ denote the coefficient of the basis element $b_\la$ in the
expansion of $f$ in this basis.   First we have a simple consequence
of Young's Rule.
\bth
\label{yr}
For any partitions $\la,\mu$ we have $[m_\mu] s_\la$ is a nonnegative
integer.  In particular,

\vs{5pt}

\eqed{
[m_\mu]s_\la=
\case{$1$}{if $\mu=\la$,}
{$0$}{if $\mu\not\le \la.$}
}
\eth

Let $\la+\mu$ denote the componentwise sum
$(\la_1+\mu_1,\la_2+\mu_2,\ldots)$.  
The next result follows from the
Littlewood-Richardson Rule and induction.
\bth
\label{lrr}
For any partitions $\la^1,\ldots,\la^r$ and $\mu$ we have
$[s_\mu] s_{\la^1}\cdots s_{\la_r}$ is a nonnegative integer.  In particular,

\vs{5pt}

\eqed{
[s_\mu]s_{\la^1}\cdots s_{\la^r}=
\case{$1$}{if $\mu=\la^1+\cdots+\la^r$,}
{$0$}{if $\mu\not\le \la^1+\cdots+\la^r$.}
}
\eth
Because of this result we call $\la^1+\cdots+\la^r$ the 
{\it dominant partition\/} for $s_{\la^1}\cdots s_{\la^r}$.

Finally, we need a result of Kirillov~\cite{kir:csg} about the product
of Schur functions, which was proved bijectively by Kleber~\cite{kle:prs} and
Fulmek and Kleber~\cite{fk:bps}.  This result can be obtained by applying the Desnanot-Jacobi Identity---also known as Dodgson's condensation formula---to the Jacobi-Trudi matrix for $s_{k^{r+1}}$.
Note that, to improve readability, we drop the sequence parentheses when a sequence appears as a subscript.

\bth[\cite{fk:bps,kir:csg,kle:prs}]
\label{kir:thm}
For positive integers $k,r$ we have

\vs{5pt}

\eqed{
\left(s_{k^r}\right)^2 - s_{(k-1)^r} s_{(k+1)^r} = s_{k^{r-1}} s_{k^{r+1}}.
}
\eth

To state our results, we need a few more definitions.  If ${b_\la}$ is
a basis for $\La(\bx)$ and $f\in\La(\bx)$ then we 
say $f$ is {\it $b_\la$-nonnegative\/} if $[b_\la]f\ge0$ for all
partitions $\la$.  Note that $m_\la$-nonnegativity is the natural
generalization to many variables of the $q$-nonnegativity definition for
$\bbR[q]$.  Also note that $s_\la$-nonnegativity implies
$m_\la$-nonnegativity by Theorem~\ref{yr}. 
\bth
The sequence $\cL^i(h_k)$ is $s_\la$-nonnegative  for $0\le  i\le 3$.  
But the sequence $\cL^4(h_k)$ is not $m_\la$-nonnegative.
\eth
\bprf
From the definition of the Schur function we have
$$
L^0(h_k)=h_k=s_k \quad\mbox{and}\quad L^1(h_k)=(h_k)^2-h_{k-1}h_{k+1}=s_{k^2}.
$$
Now Theorem~\ref{kir:thm} immediately gives
$$
L^2(h_k) = \left(s_{k^2}\right)^2 - s_{(k-1)^2}s_{(k+1)^2} = s_k s_{k^3}
$$
which is $s_\la$-nonnegative by the first part of Theorem~\ref{lrr}.
Using Theorem~\ref{kir:thm} twice gives
\begin{eqnarray*}
L^3(h_k)
&=& \left(s_k\right)^2 \left(s_{k^3}\right)^2 - s_{k-1} s_{(k-1)^3} s_{k+1} s_{(k+1)^3}\\
&=& \left(s_k\right)^2 \left(s_{k^3}\right)^2 - \left(s_k\right)^2 s_{(k-1)^3} s_{(k+1)^3}\\
& & + \left(s_k\right)^2 s_{(k-1)^3} s_{(k+1)^3} - s_{k-1} s_{(k-1)^3} s_{k+1} s_{(k+1)^3}\\
&=& \left(s_k\right)^2 s_{k^2} s_{k^4} +  s_{(k-1)^3} s_{k^2} s_{(k+1)^3}
\end{eqnarray*}
which is again $s_\la$-nonnegative.  This finishes the cases 
$0\le i\le3$.

We now assume $k\ge2$.  Computing $L^4(h_k)$ from the expression for
$L^3(h_k)$ gives the sum of the terms in the left column below.  The right column gives the dominant partition for each term, as determined by Theorem~\ref{lrr}.
\[
\begin{array}{ll}
+(s_k)^4 (s_{k^2})^2  (s_{k^4})^2 
& (8k, 4k, 2k, 2k) \\
+2 (s_k)^2 (s_{k^2}) ^2  s_{k^4} s_{(k-1)^3} s_{(k+1)^3} 
& (7k, 5k, 3k, k) \\
+ (s_{(k-1)^3})^2 (s_{k^2})^2  (s_{(k+1)^3})^2
& (6k, 6k, 4k) \\
-  (s_{k-1})^2 s_{(k-1)^2} s_{(k-1)^4} (s_{k+1})^2 s_{(k+1)^2} s_{(k+1)^4} 
& (8k, 4k, 2k, 2k) \\
- (s_{k-1})^2 s_{(k-1)^2} s_{(k-1)^4} s_{k^3} s_{(k+1)^2} s_{(k+2)^3} 
& (7k-1, 5k+1, 3k+1, k-1) \\
-s_{(k-2)^3} s_{(k-1)^2} s_{k^3} (s_{k+1})^2 s_{(k+1)^2} s_{(k+1)^4} 
& (7k+1, 5k-1, 3k-1, k+1)\\
-s_{(k-2)^3} s_{(k-1)^2} (s_{k^3})^2 s_{(k+1)^2} s_{(k+2)^3} 
& (6k, 6k, 4k)
\end{array}
\]

Now consider $\la=(7k+1,5k-1,3k-1,k+1)$, the dominant partition for the penultimate term above.  Observe that if $\mu$ is the dominant partition for any other
term, then $\la\not\le\mu$.  So, by the
second part of Theorem~\ref{lrr}, $s_\la$ appears in the Schur-basis
expansion for $L^4(h_k)$ with coefficient $-1$.  It then follows from the 
second part of Theorem~\ref{yr}, that
the coefficient of $m_\la$ is $-1$ as well.
\eprf

\section{Real roots and Toeplitz matrices}
\label{rrp}

We now consider two other (almost equivalent) settings where, in contrast
to the results of the previous section, Conjecture~\ref{bm:con} does
seem to generalize.  In fact, this may be the right level of
generality to find a proof.

Let $(a_k)=a_0,a_1,\ldots,a_n$ be a finite sequence of nonnegative real numbers.
It was shown by Isaac Newton that if all the roots of the polynomial
$p[a_k]\stackrel{\rm def}{=}a_0 + a_1 x + \cdots a_n x^n$ are real, then the
sequence $(a_k)$
is log-concave.  For example, since the polynomial
$(1+x)^n$ has only real roots, the $n$th row of Pascal's triangle is
log-concave.  It is natural to ask if the real-rootedness property is 
preserved by the $\cL$-operator.
The literature includes a
number of results about operations on polynomials which preserve
real-rootedness; for example,
see~\cite{bra:ltp,bre:ulp,bre:lus,pit:pbc,wag:tph,wy:prz}.

\bcon
\label{pla:con}
Let $(a_k)$ be a finite sequence of nonnegative real numbers.  If $p[a_k]$
has only real roots then the same is true of $p[\cL(a_k)]$.
\econ

This conjecture is due independently to 
Richard Stanley~[private communication].  It is also one of a number
of related conjectures made by Steve Fisk~\cite{fis:qdp}.
If true, Conjecture~\ref{pla:con}
would immediately imply the original Boros-Moll Conjecture.
As evidence for the conjecture, we have verified it by computer for a
large number of randomly chosen real-rooted polynomials.  We have also checked
that
$p[\cL_k^i\binom{n}{k}]$ has only real roots for all $i\le 10$ and
$n\le 40$.  It is interesting to note that  Boros and Moll's
polynomial $P_m(x)$ in equation~\ree{Pm} does not have real roots even
for $m=2$.  So if the corresponding sequence is infinitely log-concave
then it must be so for some other reason.

Along with the rows of Pascal's triangle, it appears that applying $\cL$ to the other finite lines we were considering in Section~\ref{cpt} also yields sequences with real-rooted
generating functions.
So we make the following conjecture
which implies the ``if'' direction of Conjecture~\ref{all:con}.
\bcon
For $0\le u<v$, the polynomial $p[\cL_m^i(\binom{n+mu}{mv})]$  has
only real roots for all $i\ge0$.  
\econ
We have verified this assertion for all $n\le24$ with $i\le10$ and
$v\le10$.  
In fact, it follows from a theorem of Yu~\cite{yu:ctc} that the conjecture holds for $i=0$ and all $0\le u<v$.  So it will suffice to prove Conjecture~\ref{pla:con} to obtain this result for all $i$.

We can obtain a matrix-theoretic perspective on problems of
real-rooted\-ness via 
the following renowned result of Aissen, Schoenberg and
Whitney~\cite{asw:gft}.  A matrix $A$ is said to be totally
nonnegative if every minor of $A$ is nonnegative.  
We can associate with any sequence $(a_k)$ a corresponding
(infinite) {\it Toeplitz matrix\/} $A = (a_{j-i})_{i,j\ge0}$.
In comparing the
next theorem to Newton's result, note that for a real-rooted
polynomial $p[a_k]$ the roots being nonpositive is
equivalent to the sequence $(a_k)$ being nonnegative.

\bth[\cite{asw:gft}]
\label{asw:th}
Let $(a_k)$ be a finite sequence of real numbers.  Then every root of
$p[a_k]$ is a nonpositive real number if and only if the 
Toeplitz matrix $(a_{j-i})_{i,j\ge0}$ is totally nonnegative. \hqed
\eth

To make a connection with the $\cL$-operator, note that
$$
a_k^2-a_{k-1}a_{k+1}= \left| \begin{array}{cc} a_k & a_{k+1} \\ a_{k-1} & a_k
\end{array} \right|,
$$
which is a minor of the Toeplitz matrix 
$A =(a_{j-i})_{i,j\ge0}$.  Call such a minor {\it adjacent}
since its entries are adjacent in $A$.    
Now, for an arbitrary infinite matrix 
$A = (a_{i,j})_{i,j\geq 0}$, let us define the infinite matrix $\cL(A)$ by 
$$
\cL(A) = \left( \left| 
\begin{array}{cc} a_{i,j}   & a_{i,j+1} \\
                  a_{i+1,j} & a_{i+1,j+1} 
\end{array} 
\right| \right)_{i,j \ge0}. 
$$
Note that if $A$ is the Toeplitz matrix of $(a_k)$ then $\cL(A)$ is
the Toeplitz matrix of $\cL(a_k)$.
Using Theorem~\ref{asw:th}, Conjecture~\ref{pla:con} can now be strengthened as
follows.

\bcon
\label{tnn:con}
For a sequence $(a_k)$ of real numbers,
if $A = (a_{j-i})_{i,j\ge0}$ is totally nonnegative then $\cL(A)$ is also
totally nonnegative.
\econ

Note that if $(a_k)$ is finite, then Conjecture~\ref{tnn:con} is equivalent to Conjecture~\ref{pla:con}.
As regards evidence for Conjecture~\ref{tnn:con}, consider an arbitrary $n$-by-$n$
matrix $A = (a_{i,j})_{i,j=1}^n$.  For finite matrices, $\cL(A)$ is
defined in the obvious way to be the $(n-1)$-by-$(n-1)$ matrix
consisting of the 2-by-2 adjacent minors of $A$.  
In~\cite[Theorem 6.5]{fhgj:ctp}, Fallat, Herman, Gekhtman, and Johnson
show that for $n\le4$, $\cL(A)$
is totally nonnegative whenever $A$ is.  However,
for $n=5$, an example from their paper can be modified to show that
if 
$$
A = 
\left(
\begin{array}{ccccc}
1 & t & 0 & 0 & 0 \\
t & t^2+1 & 2t & t^2 & 0\\
t^2 & t^3+2t & 1+4t^2 & 2t^3+t & 0\\
0 & t^2 & 2t^3+2t & t^4+2t^2+1 & t\\
0 & 0 & t^2 & t^3+t & t^2
\end{array}
\right)
$$
then $A$ is totally nonnegative for $t\ge0$,  but $\cL(A)$ is not
totally nonnegative for sufficiently large $t$ ($t\ge\sqrt{2}$ will
suffice).  We conclude that the Toeplitz structure would be important
to any affirmative answer to Conjecture~\ref{tnn:con}. 

We finish our discussion of the matrix-theoretic perspective with a
positive result similar in flavor to Conjecture~\ref{tnn:con}. 

\begin{prop}
If $A$ is a finite square matrix that is positive semidefinite, then $\cL(A)$
is also positive semidefinite.
\end{prop}

\bprf
The key idea is to construct the {\it second compound matrix}
$\mathcal{C}_2(A)$ of $A$, which is the array of all 2-by-2 minors of
$A$, arranged lexicographically according to the row and column
indices of the minors~\cite{hj:ma}.  

We claim that if $A$ is
positive semidefinite, then so is $\mathcal{C}_2(A)$.  Indeed, since
the compound operation preserves multiplication and inverses, the
eigenvalues of $\mathcal{C}_2(A)$ are equal to the eigenvalues of
$\mathcal{C}_2(J)$, where $J$ is the Jordan form of $A$.  If $J$ is
upper-triangular and has diagonal entries $\la_1, \la_2, \ldots,
\la_n$, then we see that $\mathcal{C}_2(J)$ is upper-triangular with
diagonal entries $\la_i \la_j$ for all $i<j$.  Since the $\la_i$'s are
all nonnegative, so too are the eigenvalues of $\mathcal{C}_2(J)$,
implying that $\mathcal{C}_2(A)$ is positive semidefinite.   

Finally, since $\cL(A)$ is a principal submatrix of
$\mathcal{C}_2(A)$, $\cL(A)$ is itself positive semidefinite.   
\eprf

\medskip

{\it Acknowledgements.}  We thank Bodo Lass for suggesting that we
approach Conjecture~\ref{bm:con} from the point-of-view of real roots
of polynomials.  Section~\ref{rrp} also benefited from interesting
discussions with Charles R. Johnson.

\end{document}